\documentclass[12pt]{amsart}
\newtheorem{theorem}{Theorem}[section]
\newtheorem{lemma}[theorem]{Lemma}    
\newtheorem{proposition}[theorem]{Proposition}

\newtheorem{definition}[theorem]{Definition}
\newtheorem{corollary}[theorem]{Corollary}

\renewcommand{\proof}{\underline{Proof.}\quad}

\newcommand{\Co}{\mbox{$\mathbb{C}$}}
\newcommand{\N}{\mbox{$\mathbb{N}$}}

\begin{document}
   \title[Multipliers and dual operator algebras] 
{Multipliers and dual operator algebras}

\author{David P. Blecher}
\address{Department of Mathematics\\University
 of Houston\\Houston,
TX 77204-3476 } \email{dblecher@math.uh.edu} 
 \thanks{August 2000.  Revision of January 10, 00.
Some results herein were presented at the Conference on 
Operator Function Theory, Ambleside, Sept. 8, 00.} 
\thanks{* Supported by a grant from the NSF} 

\vspace{30 mm}

\maketitle

\vspace{40 mm}

\begin{abstract} In a previous paper we showed how the 
main theorems characterizing operator algebras and 
operator modules, fit neatly into the framework of the
`noncommutative Shilov boundary', and more particularly
via the left multiplier operator algebra of an 
operator space.   As well as giving new characterization
theorems, the approach  of that paper
allowed many of the hypotheses of the earlier  
theorems to be eliminated.
Recent progress of the author with Effros 
and Zarikian now enables weak*-versions of these
characterization theorems.  For example, we prove 
a result analogous to Sakai's famous
characterization of 
von Neumann algebras as the $C^*$-algebras with predual,
namely, that the $\sigma$-weakly closed unital 
(not-necessarily-selfadjoint) subalgebras of $B(H)$
for a Hilbert space $H$, are exactly the unital operator 
algebras which possess an operator space predual.
This removes one of the hypotheses from an earlier 
characterization due to Le Merdy.
We also show that the multiplier operator algebras
of dual operator spaces are dual operator algebras.
Using this we refine several known results 
characterizing dual operator modules.
\end{abstract}

\pagebreak
\newpage

\section{Introduction.}  

This paper is an application of 
non-commutative M-structure.
Classically, M-structure is a theory which grew out 
of the seminal paper \cite{AE} (see \cite{HWW} for more details).  
This theory has many facets, one of
which is the theory and application of M-ideals (and, by 
duality, L-ideals).  
Recall that an arbitrary Banach space $X$ can be
embedded linearly isometrically in an `M-space' $C(K)$ (in many 
ways).  There are also constructions that are `M in nature' 
such as the $\oplus^\infty$ direct sum.
Part of classical `M-theory' is concerned with what 
traces of the structure of the $C(K)$ spaces, or 
the $\oplus^\infty$ direct sum, are visible in a general
$X$.   Via such considerations,
a simple looking metric condition in $X$ will often force a
powerful conclusion of a structural nature.
 Another aspect of classical M-structure is the multiplier
and centralizer algebras of a Banach space $X$.  These are, respectively,
a function algebra $Mult(X)$, and commutative C*-algebra $Cent(X)$,
which may be associated with any Banach space $X$, and which, in 
some sense, control the M-structure of $X$.  It is these algebras
 that we concentrate on here.  In an earlier
paper \cite{B} we defined noncommutative generalizations of
$Mult(X)$ and $Cent(X)$, 
which will be explained in more detail presently.

The main application of this paper is to theorems characterizing
basic algebraic structures that are of interest to operator 
algebraists.  A fundamental theorem in the field of
operator algebras is the Gelfand and Naimark characterization of 
$C^*$-algebras (see \cite{Ped} for example).  
For nonselfadjoint operator algebras the 
appropriate characterization is the BRS theorem \cite{BRS}
which may be stated as follows: the 
norm closed subalgebras of $B(H)$ containing $I_H$,
are exactly the unital algebras $A$ which are also operator spaces
such that the identity in the algebra has norm 1, and 
also the condition:
$$\Vert a b \Vert_n \leq  \Vert a  \Vert_n \Vert b \Vert_n $$
holds, for matrices $a, b \in M_n(A)$.   Here the words ``are exactly''
mean ``up to a homomorphism (i.e. an isomorphism)
which is completely isometric''.   
We will give a short 
proof of BRS, and indeed of a new 1-sided version
of BRS, at the end of this introduction.

With this theorem in mind, when we refer henceforth to an
operator algebra in this paper, we mean either 1) a `concrete 
operator algebra', i.e. a norm closed subalgebra
of $B(H)$ with an identity of norm 1, or 2) an operator space 
which is a unital algebra completely isometrically isomorphic
to a concrete operator algebra.  Note that any norm closed 
subalgebra $A$
 of $B(H)$ with an identity of norm 1, may be assumed to have 
identity $I_H$, by considering the completely isometric 
homomorphism on $A$ which restricts $A$ to act on the 
Hilbert space $(AH)^{\bar{}}$.  In this paper we will not
consider operator algebras with contractive approximate 
identities, since they 
are automatically unital in 
the w*-situation (see the end of \S 2).  It is worth
remarking that all function algebras (i.e. unital
norm closed subalgebras of $C(K)$) are operator algebras.
Indeed these may be characterized as the  operator algebras
which are also $MIN$ operator spaces \cite{Bcomm}.

There are also known characterization theorems similar to 
`BRS' appropriate to the class of
`operator modules' (which are defined below)
 up to completely isometric module 
isomorphism (see \cite{CES} for the
earliest such result).   

In \cite{B} we were able to fit these
 characterization theorems into a coherent context 
using the multiplier algebras just mentioned.  This approach
had three advantages of 1) being a more unified treatment of
 these theorems,  2) allowing many of the hypotheses of 
the earlier theorems to be eliminated, and 3) giving new results.
However it was not at all clear from \cite{B} whether the new
approach could be applied to the weak*-versions of such 
characterization theorems.  For example, for $C^*$-algebras
the appropriate weak*-version is Sakai's famous
theorem, which 
characterizes von Neumann algebras as exactly the 
$C^*$-algebras with predual.   Here the words ``as exactly''
mean ``up to a *-homomorphism which is also a
w*-w*-homeomorphism''.  Sakai's result is of course 
a noncommutative generalization of the classical 
characterization of $L^\infty$ spaces as exactly the $C(K)$
spaces with predual.  One of the main results of the 
present paper is the following version of Sakai's 
theorem for (not-necessarily-selfadjoint) operator algebras.

\begin{theorem}  \label{Blm} Let  $A$ be an operator
algebra with an identity of norm 1, which is also
a dual operator space (i.e. $A = X^*$ completely 
isometrically for an
operator space $X$).  Then $A$ is w*-w*-homeomorphically and  
completely isometrically isomorphic  to a
$\sigma$-weakly closed (i.e. w*-closed) subalgebra
$B$ of  $B(H)$ for some Hilbert space $H$.   Conversely,
any $\sigma$-weakly closed subalgebra (or indeed subspace)
$B$ of  $B(H)$ is a dual operator space.
\end{theorem}

Of course the isomorphism here is also multiplicative,
linear, and w*-w*-homeomorphic.   The converse direction of
the theorem is an old and quite simple
result (see 4.2.2 in \cite{ERbook}).   
          
By virtue of the theorem, we shall henceforth refer to a
$\sigma$-weakly closed subalgebra of  $B(H)$ containing 
an identity of norm 1 as a 
`concrete dual operator algebra', and an operator 
algebra which has an operator space predual as
a `dual operator algebra'.

Some important remarks
are in order.  Firstly, one may assume as before
 that the identity of norm 1 in a concrete 
dual operator algebra
is indeed $I_H$, by restricting 
$A$ to $(AH)^{\bar{}}$.  This restriction is clearly a
completely isometric, weak*-continuous, homomorphism,
and  hence it is
a w*-w*-homeomorphism (see Lemma \ref{KS} (3)).

The second important remark is that it turns out that the 
hypothesis that $A$ be a dual operator space, as opposed 
to a dual Banach space, is necessary for the conclusion
of the theorem
to hold.  We will discuss an example
in \S 2.

The third important remark is that Theorem \ref{Blm} is
a refinement of the following theorem of Christian Le 
Merdy, which appeared first in \cite{LM3}.
A very
natural and elementary proof of it appears in \cite{LM4}
(Proposition 3.4 and Remark 3.5).

\begin{theorem}  \label{LM} (Le Merdy) Let $A$ be as in the statement
of \ref{Blm}, but with the additional hypothesis that 
the multiplication on $A$ is separately weak*-continuous.
Then  $A$ is w*-w*-homeomorphically 
completely isometrically isomorphic  to a concrete dual 
operator algebra.
\end{theorem}

Thus our proof of \ref{Blm} will consist of showing that the 
multiplication on an operator algebra which is a dual 
operator space, is automatically 
separately weak*-continuous.
 
We shall prove Theorem \ref{Blm} in \S 2.  In \S 3 we will use this
to deduce that the left multiplier algebra of any dual 
operator space is a dual operator algebra.   We will 
also give a quick `M-structure' proof of an operator space
variant on important results of Zettl, and Effros, Ozawa and 
Ruan \cite{Zettl,EOR}.    
In \S 4 we extend several results from 
\cite{B} to the weak* situation.  
Before we describe some such results, we remark
that hitherto this was quite problematic, because 
the `Hamana theory' or  `Shilov boundary' approach 
did not seem to work well for dual spaces.  In particular
this is because 
the noncommutative Shilov boundary ${\mathcal T}(X)$ and 
injective envelope $I(X)$ of a dual operator space $X$ 
did not seem related to von Neumann algebras in any 
sense we were able to use.  However
the various left multiplier operator algebras of $X$ {\em are}
dual algebras (as we show here in \S 3 and in \cite{BEZ}).

To describe the work in the last section we recall the following

\begin{definition} \label{omd}
An {\em operator module} is an operator space $X$, together with a
left module action $A \times X \rightarrow X$ which satisfies:
\begin{itemize}
\item [(i)]
$e x = x$  for all $x \in X$,
\item [(ii)]  $\Vert a x \Vert_n \leq \Vert a \Vert_n \Vert x
\Vert_n$  for all matrices $a \in M_n(A), x \in M_n(X)$ of any
size ($n < \infty$).
\end{itemize}
\end{definition}

Here $e$ is the identity of $A$.  We have not said anything about
what $A$ is, traditionally it is a unital $C^*$-algebra or
operator algebra, but we will want to generalize this in \S 4.
Note that (ii) may be rephrased as saying
that the module action is completely contractive as a bilinear map
(in the sense of \cite{CS,PS}).

As motivation for our work in  \S 4 we recall the following

\begin{theorem} \label{ern}  (Effros-Ruan \cite{ERbimod})  Let
$M$ be a $W^*$-algebra, and $X$ a left operator $M$-module which
is also a dual operator space, such that the module
action $M \times X \rightarrow X$ is separately w*-continuous.
Then there exist Hilbert spaces $H, K$, a unital normal
*-representation $\pi : M \rightarrow B(K)$, and a
weak*-homeomorphic complete isometry $\Phi : X \rightarrow
B(H,K)$, such that $\Phi(bx) = \pi(b) \Phi(x)$
for all $x \in X, b \in M$.
\end{theorem}

A slightly more general result with a different proof 
may be found in \cite{BLM}. 
Effros and Ruan call an $X$ satisfying the hypotheses of this result,
a {\em normal dual left operator
module} over $M$.
A little thought  will convince the reader that
this result is an abstract characterization of
w*-closed subspaces of $B(H,K)$ which are left invariant
under multiplication by a von Neumann algebra acting
on $K$.

We are able in \S 4 to improve on this useful
characterization in four ways:
1) we show that the hypothesis that the action is separately w*-continuous
it not necessary, one only needs that the action be
w*-continuous in the first variable and then
w*-continuity in the second variable is automatic;
2) we show that it can be
arranged so that $H, K$ and
 $\Phi$ only depend on $X$, and not on $M$ or
the particular action; 3) in the bimodule version of
this theorem our general approach
(\cite{B} \S 5) shows that the left and right module
actions on $X$ {\em automatically} commute with
each other; and 4) we will generalize this result
to allow $M$ and its action on $X$ to be replaced by any
dual operator space  which is a unital algebra
such as perhaps a group algebra.  Indeed $M$ need not even
be an algebra, we will be able to characterize `oplications'
$m : Y \times X \rightarrow X$
which are w*-continuous in appropriate variables.
An `oplication' (see \cite{B} \S 5)
is a bilinear map $Y \times X \rightarrow X$,
for operator spaces $X, Y$, which
also satisfies (i) and (ii) of definition \ref{omd},
where $e$ is a fixed element of $Ball(Y)$.
We say that $e$ is the `unit' of $Y$.

Thus we have 
extended the main results of \S 5 of \cite{B}  
to dual modules and oplications.

We thank E. G. Effros, Christian Le Merdy, and N. Weaver 
for helpful comments. 
 
We end this introduction with some notation and basic 
results.  We first list a well known
and  basic functional analysis
fact which we will use very frequently:

\begin{lemma} \label{KS} (Krein-Smulian Theorem)  \begin{itemize}
\item [(1).]  Let 
$X$ be a dual Banach space, and $Y$ be
a linear subspace of $X$.  Then 
$Y$ is w*-closed in $X$ iff  $Ball(Y)$ is closed in
the  w*-topology on $X$.  In this case $Y$ is also a
dual Banach space, with predual $X_*/ Y_\perp$, and the  
inclusion of $Y$ in $X$ is w*-continuous.
\item [(2).]  A linear bounded map $T$ between
dual Banach spaces is
weak*-continuous if and only if whenever $x_i \rightarrow x$
is a bounded net converging
weak*- in the domain space, then
$T(x_i) \rightarrow T(x)$ weak*.   
\item [(3).]   Let $X$ and $Y$ be dual Banach spaces, and 
$T : X \rightarrow Y$  a w*-continuous linear isometry.
Then $T$ has w*-closed range $V$ say, and $T$ is a
w*-w*-homeomorphism onto $V$.
\end{itemize}  
\end{lemma} 

\begin{proof}  (1) and (2) may be found in any book on  basic 
functional analysis.  (2) is often stated for functionals $\phi$
but the result as stated follows from this by considering 
$\phi \circ T$.  Item (3) is found in fewer books,
the proof is quite obvious from (1).  For by 
(1), $V$ is clearly w*-closed in $Y$, and the restriction of
$T$ to $Ball(X)$ then takes w*-closed (and thus 
w*-compact) sets to w*-compact (and thus w*-closed)
sets in $V$.  Thus the inverse of $T$  
restricted to the ball is  w*-continuous,
so $T^{-1}$ is w*-continuous by (2).
\end{proof}

\vspace{5 mm}

We will use 
the term $W^*$-algebra for a $C^*$-algebra with predual.
In view of the aforementioned theorem of Sakai this
is `the same as' a von Neumann algebra.

We now turn to operator spaces and left multipliers.  An operator 
space  is a linear subspace of $B(H)$ for a Hilbert space $H$.
In this paper all of our operator spaces are norm complete.
Equivalently, there is an abstract characterization due to 
Ruan.  We refer the reader to the books \cite{ERbook,Pis}
for more details on operator spaces.
Any details needed which are omitted from
 these books may be found in \cite{BP,Bsd}.
Details on completely bounded and completely positive 
maps may be found in
\cite{P}.  We write $M_{nm}(X)$ for the operator space of 
$m \times n$ matrices with entries in $X$, and 
$C_n(X) = M_{n,1}(X)$ and $R_n(X) = M_{1,n}(X)$ as usual. 
For a map $T : X \rightarrow Y$ we write
$T^{(n)}$ for the associated amplification 
$M_n(X) \rightarrow M_n(Y)$, that is $T \otimes I_n$.  
The `completely' prefix to a property means that it
passes to every $T^{(n)}$, thus for example
$T$ is completely isometric means that each $T^{(n)}$ is
isometric.   

We now turn to duality; and here we have a notational 
problem, the use of the symbol $*$, which is used 
for three different things in this paper.  Namely,
we have the dual space $X^*$ of a space $X$; the  
adjoint or involution $S^*$ of an operator on a Hilbert space, and 
the adjoint operator $R^* : Y^* \rightarrow Z^*$ of 
an operator $R : Z \rightarrow Y$.   We are forced by reasons
of personal
taste to leave it to the reader to determine which is 
meant in any given formula; although the third of these will
only occur once or twice.  We recall that if $X$ is
an operator space then so is $X^*$; its matrix norms
come from the identification $M_n(X^*) \cong CB(X,M_n)$.
A {\em dual operator space} is one which is completely 
isometrically isomorphic to the operator space dual of 
another operator space.   Any $\sigma$-weakly closed 
subspace of $B(H)$ is a dual operator space with 
predual determined by the predual of $B(H)$.
Conversely, any dual operator space 
is linearly completely
isometrically w*-w*-homeomorphic to a
$\sigma$-weakly closed subspace of some $B(H)$.
If $X$ is
an operator 
space which is the dual Banach space of $Y$, then there is 
 at most one operator space structure (i.e. matrix norms) on
$Y$ with respect to which it is possible for 
$X = Y^*$ completely isometrically (since $Y \subset Y^{**} = 
X^*$).  But one should 
be warned: there may in fact be no such operator space structure
on $Y$  (see \cite{LM1,EOR} 
for examples).
However there is always one if, further,
 $X$ is a C*-algebra \cite{Bsd}
(but not necessarily if $X$ is a  
nonselfadjoint operator algebra) or a $MIN$ or $MAX$ space.

It will be important for us 
to note that if $X$ is a  dual operator space
then so is $M_n(X)$, with predual equal to the 
operator space projective tensor product of the 
predual $T_n$ of $M_n$, and the predual of $X$.  The 
duality pairing may be taken to be the obvious 
one.  From this
it follows that:

\begin{lemma}  \label{mac}  If $X$ is a  dual operator space,
and $x_i$ is a net in $M_n(X)$, then $x_i \rightarrow
x \in M_n(X)$ weak*- in $M_n(X)$, if and only if
each entry in $x_i$ converges weak*- in $X$ to the 
corresponding entry in $x$.
\end{lemma}

We also recall \cite{BP,ERbook}
that for a dual operator space $X = Y^*$, the
space $CB(X)$ is canonically a dual operator space, with
predual the operator space projective tensor
product of $X$ and $Y$.  A bounded net $T_i \in CB(X)$
thus converges in the w*-topology to $T \in CB(X)$ if and
only if $T_i(x) \rightarrow T(x)$ w*- in $X$, for every
$x \in X$.  From this it is obvious that the multiplication
on $CB(X)$, viewed as a bilinear map,
 is weak*-continuous in the first variable.    

We will not use much about operator systems, but refer 
the reader to \cite{ERbook,P} for details.   For the 
purposes of this paper we will define a {\em dual
operator system} to simply be an operator system which is
a  dual operator space.   A unital operator space is 
a pair $(X,e)$, consisting of an operator space $X$ which 
may be completely isometrically linearly embedded in a 
unital $C^*$-algebra $B$, with $e$ embedded as $1_B$.  We
call $e$ the `unit' or `identity' of $X$, and identify
unital operator spaces up to `unital complete isometry'.
Define a subset  $\Delta(X) = 
\{ b \in X : b^* \in X \}$ of $B$.   Then $\Delta(X)$ 
is an operator system which is
independent (up to unital complete order isomorphism)
of the actual $B$ containing $X$ unitally.
This follows from Arveson's result \cite{Arv1}
that the space $\{x + y^* : x,y \in X \}$
 is well defined up to complete
order isomorphism.  Therefore given a unital operator space
$X$, the notation $\Delta(X)$ makes sense even if no
particular $B$ containing $X$ is specified.

Notice further, that if $X$ in the last paragraph
is an operator algebra $A$ with identity of norm
$1$,  then we may suppose that $B$ above is a 
$C^*$-algebra containing $A$ as a unital subalgebra.
Hence the space $\Delta(A)$ above is a closed
$C^*$-subalgebra of $B$.  This $C^*$-algebra 
$\Delta(A)$ is independent up to *-isomorphism of the 
particular  $B$.
This is well-known fact.       
 
We turn to the space $M_l(X)$ of left multipliers, and the 
space $A_l(X)$ of left adjointable maps, on 
an operator space $X$.  These were defined first in
\cite{B}, although earlier authors had considered 
variants valid for operator systems\footnote{There appears to
be no way to define operator space left multipliers
via left multipliers of an associated 
operator system.}\cite{KW,Ki}.  
They are natural generalizations of classical spaces, for example
$M_l(X)$ generalizes the space of  multipliers of a 
Banach space\footnote{These
may be defined to be the linear maps $T$ on the Banach
space $X$ of the form $Tx = gx$ for $x \in X$, where $g \in C(K)$
for some $C(K)$  containing $X$ isometrically.} \cite{AE,HWW}.

The reader not willing to refer to \cite{B} for 
the original definitions, or to \cite{B,BEZ,BPnew} for 
various equivalent formulations, may take the definitions from 
the following result from \cite{B}.  

\begin{theorem}  \label{chmu}  
Let $X$ be an 
operator space, and $T : X \rightarrow X$ a linear map.  Then:
\begin{itemize} 
\item [(1).]  $T \in M_l(X)$ if and only if 
there exists a Hilbert space $H$, an $S \in B(H)$, and a
completely isometric
linear embedding $\sigma : X \rightarrow B(H)$ such that
$\sigma(Tx) = S \sigma(x)$ for all $x \in X$. 
\item [(2).]  $T \in A_l(X)$ if and only if there
exist $H, S , \sigma$ satisfying all the conditions of (1),
and also $S^* \sigma(X) \subset \sigma(X)$.
\end{itemize}
We have that $M_l(X)$ is  an operator
algebra, with algebra structure inherited from $B(X)$
and with 
norm the least value of 
$\Vert S \Vert$ possible in (1).  This least value is
achieved.
This operator algebra has identity which is the
identity operator on $X$.  Also,
 $A_l(X) =  \Delta(M_l(X))$, 
 and this is a $C^*$-algebra.
     \end{theorem} 

\vspace{3 mm}

In (1) we can replace $B(H)$ by $B(H,K)$ with no loss, or 
we can replace $B(H)$ with a $C^*$-algebra $A$,
 or a $C^*$-module.   
 
\begin{corollary} \label{her}  If $T$ is a left multiplier 
of an operator space $X$, then the `multiplier norm' of 
$T$ is greater than, or equal to, its completely bounded norm.
Also, if $Y$ is a closed subspace 
of $X$ with $T(Y) \subset Y$, then $T_{|_Y}$ is a
left multiplier of $Y$, with a smaller (or equal) `multiplier
norm' than that of $T$.    If in addition,
$T$ is adjointable on $X$ and $T^*(Y)  \subset Y$, then 
$T_{|_Y}$ is adjointable on $Y$.   
\end{corollary}

We will therefore consider $M_l(X)$ as a unital 
subalgebra of $CB(X)$, but 
possessing a possibly
larger `norm'.  On $A_l(X)$ the multiplier norm 
coincides with the `cb-norm', as is shown in \cite{B}
\S 4.  We 
assign matrix norms to $M_l(X)$ via the relation  
$M_n(M_l(X)) \cong M_l(C_n(X)) \cong M_l(M_n(X))$, and 
similarly for $A_l$.    
There are many ways to see this last relation,
for example using the proof of 4.8 in
\cite{B}, or by
using the `$I_{11}$' definition of $M_l$ in \cite{BPnew}.

One alternate characterization of $M_l(X)$ from 
\cite{B} shows that 
for unital operator spaces $X$, we have 
$$M_l(X) \cong \{ a \in C^*_e(X) : a J(X) \subset J(X) \} $$
where $C^*_e(X)$ is the $C^*$-envelope \cite{Ham1}, and
$J$ is the (unital)
embedding of $X$ inside $C^*_e(X)$.  
The isomorphism here is the map taking $a$ in the right hand
set to the map $x \mapsto J^{-1}(aJ(x))$ on $X$.  This map
clearly has inverse which is the map $T \mapsto J(T(1))$.
Since $1 \in X$
it follows that $M_l(X) \subset X$ in this case.  More 
specifically, the map $M_l(X) \rightarrow X$ given 
by $T \mapsto T(1)$ is a unital complete isometry.
Since $\Vert T \Vert_{cb} \geq \Vert T(1) \Vert$, we see
that for unital operator spaces
 the canonical completely contractive
inclusion $M_l(X) \rightarrow 
CB(X)$ is completely isometric.  In particular
for a unital operator algebra $A$, we have 
$M_l(A) \cong A$ completely isometrically isomorphically.  Of 
course viewed as a subset of $CB(A)$, these left
multipliers
are exactly the maps $T(b) = ab$ for fixed $a \in A$.

The following, our most recent
characterization of left multipliers, is also the deepest
and no doubt most useful.  It
is a generalization of a result of 
W. Werner characterizing multipliers on
operator systems \cite{WW}.

\begin{theorem} \label{BEZr} \cite{BEZ}
 A linear map $T : X \rightarrow X$
on an operator space is a  left multiplier of multiplier
norm $\leq 1$, if and only if
the following  map $\tau_T$ is completely contractive
on $C_2(X)$: 
$$\left[ \begin{array}{c} x \\ y \end{array}
\right] \mapsto \left[ \begin{array}{c} Tx \\ y \end{array}
\right] $$  
\end{theorem}

Note this last theorem immediately gives a proof of the `BRS'
theorem mentioned in the second 
paragraph of our paper.  More generally
it gives  a proof of the oplication theorem from \cite{B}
(from which BRS follows in one line as we pointed out in \cite{B}
5.6).

\begin{theorem} \label{oth}  (The oplication theorem).
Let $m : Y \times X \rightarrow X$ be an oplication.
Then there exists a (necessarily unique) linear complete 
contraction $\theta : Y \rightarrow M_l(X)$ such that 
$\theta(e) = I$ and 
$\theta(y)(x) = m(y,x)$ for all $y \in Y, x \in X$.
If in addition $(Y,e)$ is an operator system then 
$\theta$ maps into $A_l(X)$, whereas if 
$(Y,e)$ is a unital algebra then 
$\theta$ is a homomorphism if and only if 
$m$ is a module action.  \end{theorem}

\begin{proof}
If $m : Y \times X \rightarrow X$ is an oplication,
consider 
the associated linear map $\tilde{m} : Y \rightarrow B(X)$.
Note that 
$$\tau_{\tilde{m}(y)} \left( \left[ \begin{array}{c} x \\
x' \end{array} \right] \right) = \left[ 
\begin{array}{c} m(y,x) \\ x' 
\end{array} \right] \; ,$$
for $x, x' \in X$.  The last matrix may be regarded as 
the formal product, via $m$, of the diagonal matrix $y \oplus e$
and the column in $C_2(X)$ with entries $x$ and $x'$.
By the hypothesis, one clearly sees that $\tau_{\tilde{m}(y)}$
is contractive, and an analogous argument with matrices shows that 
it is completely contractive.  Thus by \ref{BEZr} we see that
$\tilde{m}$  is a contractive map into
$M_l(X)$. Given 
$[y_{ij}] \in Ball(M_n(Y))$, we consider $[\tilde{m}(y_{ij})]$ as 
a map $S : C_n(X) \rightarrow C_n(X)$, and argue as above using 
matrices to see that $S \in Ball(M_l( C_n(X)))$.  Thus it follows 
that $\theta : y \mapsto \tilde{m}(y)$ is 
a complete contraction from $Y 
\rightarrow M_l(X)$.  Alternatively the complete contractivity
may be deduced from the contractive case, and the 
relation $M_l(M_n(X)) \cong M_n(M_l(X))
\cong M_l(C_n(X))$.

It is clear that $\theta(y)(x) = m(y,x)$ for
all $y \in Y, x \in X$.         
The last few statements are immediate (see 5.2 in \cite{B}).
\end{proof}  

\vspace{3 mm}
 
We remark that
Paulsen has recently found a much simpler proof of 
\ref{BEZr}.  Combining this with the last proof
gives a short, and perhaps the best, 
route to the oplication and BRS theorems. 
Indeed this route gives a very interesting extension 
of BRS to operator algebras with a 1-sided identity:
  
\begin{theorem}  \label{1sBRS}  Let $A$ be an operator space
which is an algebra with a right identity of norm 1 
(or right c.a.i.).  Then $A$
is completely isometrically isomorphic to a concrete operator 
algebra, if and only if we have 
$$ \Vert (x \oplus Id_m) y \Vert_{n+m} \leq
\Vert x \Vert_{n} \Vert y \Vert_{n+m} $$
for all $n, m \in \N$ and 
$x \in M_n(X), y \in M_{n+m}(X)$.
\end{theorem}

To explain the notation of the theorem, we have written
$Id$ for a formal identity,
thus the expression $(x \oplus Id_m) y$ below means that the
upper $n \times (m+n)$-submatrix of $y$ is left 
multiplied by $x$, and the lower submatrix is left alone.

\begin{proof}  Proceed exactly as in the proof above to obtain a
complete contraction $\theta : A \rightarrow M_l(A)$ 
such that $\theta(a)(b) = ab$.   One does need to permute the entries
of the matrices involved in the calculation as to get them
into the form of the hypothesis, but this is a simple matter.
The result  may then be
completed as in our 1-line proof of BRS  from \cite{B} \S 5:
for example $\Vert \theta(a) \Vert
\geq \Vert a 1 \Vert = \Vert a \Vert$.
\end{proof}   

\vspace{3 mm}

Clearly the same idea gives a
version of 
`oplication' theorem which requires no `identity element'
$e$, which needs a similar hypothesis to that of \ref{1sBRS}.

Finally, in connection with Theorem \ref{BEZr} we also point out 
again here that it gives a method to recover the product in a 
unital operator algebra
from the linear operator space structure alone.   What is
more important, it appears to be a simple enough method to
be useful in practice; for example it works well in conjunction
with weak*-topologies, etc.   Namely, suppose
that $A$ is an operator algebra with an identity of norm $1$, but 
that we have forgotten the product on $A$.  Let us assume for a
moment that we do remember the identity element $e$. 
 Form $M_l(A)$ using \ref{BEZr}  (or if $A$ is a 
$C^*$-algebra, using Lemma 4.5 in \cite{BEZ},
(which is valid for $C^*$-algebras by going to
 the second dual)),
and define $\theta : M_l(A) \rightarrow A$ by 
$\theta(T) = T(e)$.  Then the product on $A$ is 
$ab = \theta(\theta^{-1}(a) \theta^{-1}(b))$.

If we have also forgotten the specific identity element $e$, then 
one may only retrieve the product on $A$ up to a unitary $u$ with 
$u, u^{-1} \in A$.   Such unitaries form a group.
Indeed such unitaries may be characterized 
by the Banach-Stone theorem for operator algebras (see e.g. the 
last page of \cite{B}, or 
\cite{Kad} for the $C^*$-algebra case) as 
the elements $x_0$ with the 
property that the map $\pi : T \mapsto T(x_0)$ is
a completely isometric surjection 
$M_l(A) \rightarrow A$.  
(Indeed if $A$ is a $C^*$-algebra one only needs this to
be an isometry, by Kadison's result
\cite{Kad}.  We remark that by 
Lemma 4.5 in \cite{BEZ} the unitaries in a $C^*$-algebra
correspond 
to linear $T : A \rightarrow A$ such that 
$\tau_T$ is a surjective isometry.
However in this case there are other Banach space
characterizations of unitaries - C. Akemann has shown me one 
such).   Given such an $x_0$ and $\pi$, we may again 
recover the product as
$ab = \pi(\pi^{-1}(a)\pi^{-1}(b))$.  This is the operator
algebra product on $A$ which has this unitary 
as the identity.  This is all fairly easy to see from the 
Banach-Stone theorem and our notion of the left multiplier
algebra.  Nonetheless it is quite interesting that we can
recover the product in this way.
                    
\section{A characterization of dual operator algebras}

In this section $A$ is an operator algebra with identity of 
norm 1 possessing an operator space predual $A_*$ say.
We need to show that: 

\begin{theorem}  \label{swc} If $A$ is a unital operator algebra
with operator space predual, then the 
multiplication on                         
$A$ is separately w*-continuous.  \end{theorem}              

We will consider the completely isometric unital
homomorphism $L :  A \rightarrow CB(A)$ given by 
$L(a)(b) = ab$ for $a, b \in A$.  As  we said in \S 1,
the range of $L$ is
$M_l(A)$, and quite clearly in this case,
$M_l(A) \subset CB(A)$ completely isometrically.
As we also said in \S 1,
we know that  $CB(A)$ is the dual of
the projective tensor product of 
$A$ and $A_*$, completely
isometrically.  Thus a bounded net
$T_i \rightarrow T$ weak*- in
 $CB(A)$ if and only of $T_i(b) \rightarrow T(b)$ w*- in 
$A$ for each $b \in A$.  The multiplication 
on $CB(A)$ is weak*-continuous in the first variable, quite
obviously, and we will show that this property descends to 
its subalgebra $A$.  

\begin{lemma} \label{wcl}  $M_l(A)$ is weak*-closed in
$CB(A)$.   Therefore $M_l(A)$ has a predual 
which gives $M_l(A)$ a weak*-topology coinciding 
with the relative w*-topology
on $M_l(A)$ inherited from $CB(A)$.
\end{lemma}

\begin{proof}   We will use Krein-Smulian.   
Let $T_i \in Ball(M_l(A))$ be  
a net
converging
weak* in $CB(A)$ to $T \in Ball(CB(A))$, say.   
Using Theorem \ref{BEZr}
we will now check that $T$ is  in 
$Ball(M_l(A))$.  For suppose that $v \in Ball(C_2(A))$.  Write $v$ 
as a column $[x \; \;  y]^t$. 
Then $w_i = [ T_i(x) \; \;  y]^t \in Ball(C_2(A))$
by \ref{BEZr}.   Let $w = [ T(x) \; \;  y]^t$.
Consider
the operator space predual $Z$ of $C_2(A)$; 
 any $G \in Ball(Z)$ is 
given by a pair $[\phi \; \; \psi]$ of functionals in 
$A_*$.  The duality pairing is:
$$ \langle G , v \rangle \; =  \langle \phi , x \rangle \; 
+ \langle \psi , y  \rangle \; \; . $$
 Since 
$$|\langle G , w_i \rangle | = |\langle \phi , T_i(x)
\rangle \;
+ \langle \psi , y  \rangle | \; \leq \; 1 \; \;  $$
 in the limit we have that  
$$|\langle  G  ,   w \rangle|
= |\langle \phi , T(x)
\rangle \;
+ \langle \psi , y  \rangle | \; \leq \; 1 \; .
$$ 
Thus $\tau_T$ is contractive.  A similar argument, picking
$G$ in the unit ball of the predual of 
$M_{2n,n}(A)$, and using
Lemma \ref{mac},  shows that 
 $\tau_T$ is completely  contractive.  Thus
$T \in Ball(M_l(X))$, so that  $Ball(M_l(X))$ is w*-closed.
\end{proof}

\vspace{5 mm}

\begin{proof} (Of Theorem \ref{swc}:)
By symmetry it is enough to show that
the multiplication  is weak*-continuous in the first variable.
One can see this symmetry
by considering the `opposite algebra' $A^{op}$ (c.f.
\cite{Bcomm} Proposition 1) which has
operator space predual which is $A_*$ with the transposed matrix
norm structure.
                 
From the lemma it is evident that 
the map $L^{-1} : M_l(A) \rightarrow A$ is weak*-continuous 
with respect to the weak*-topology of the Lemma.  For,
if $L(a_\lambda) \rightarrow L(a)$ w*- in $M_l(A)$,
then $
a_\lambda = L(a_\lambda)(1)  \rightarrow L(a)(1)
= a$ w*- in $A$.
By Lemma \ref{KS} (3), $L$ is 
 weak*-continuous.  But this says exactly that 
the multiplication on $A$ is  weak*-continuous 
in the first variable.
\end{proof}

\vspace{5 mm}
 
We next show
that the hypothesis that $A$ is a dual operator space
is necessary to get the full conclusion
of Theorem \ref{Blm} in general.  Note however
that if
$A$ is a function algebra, i.e. a $MIN$ space (see 
\cite{Bcomm}), then 
if $A$ has a Banach space predual then it has an operator 
space predual \cite{Bsd}.  In fact it 
is quite tricky to 
come up with examples of Banach space preduals of 
an operator space, which are not operator 
space preduals.  

Another important remark is that, in contrast to 
Sakai's theorem, the preduals of 
nonselfadjoint operator algebras are not necessarily 
unique.  See \cite{Rupr} for more detailed information.
As pointed out by Derek Westwood, this nonuniqueness may
be seen quickly 
from the fact that Banach space preduals are 
not unique, together with the well known Arveson $2 \times 2$
matrix trick (played also in the proof below).  Nonetheless,
the automatic separate w*-continuity of the multiplication 
makes some kind of uniqueness statement, which perhaps
may be more fully exploited.    

\begin{proposition} \label{lme}  (cf. \cite{LM3}) There exists an
operator
algebra $A$ with identity of norm 1, i.e. a closed subalgebra
of some $B(H)$ containing $I_H$, which is the dual of
a Banach space $X$,
and for which the multiplication on $A$ is separately w*-continuous,
but for which there is no
weak*-homeomorphic completely isometric isomorphism of $A$
onto  a dual operator algebra.  Indeed there exists 
an operator
algebra $A$ with identity of norm 1, with no operator space
predual at all, but which is a dual Banach space.
\end{proposition}

\begin{proof}
Take any operator space $Y \subset B(H)$ which is the dual of
a Banach space $X$, such that $Y$ with the associated w*-topology 
derived from $X$ is 
not completely isometrically
isomorphic  via a w*-w*-homeomorphic linear
map to a dual operator space
(see \cite{LM1} for such an example, and \cite{EOR}
for an example which has a unique Banach space predual).  
Build 
the `canonical $2 \times 2$ unital
operator algebra' $A \subset B(H \oplus H)$ 
which has $Y$ contained 
completely isometrically as the 1-2-corner and scalars
on the diagonal.   This is an adaption of an example from
\cite{LM3} where Le Merdy puts zeroes on the diagonal.
Suppose that
 $i : Y \rightarrow \ell^\infty(I)$ is
any isometric w*-homeomorphic linear embedding (such embeddings
exist since any Banach space is a quotient of a $\ell^1(I)$).
Then $A$ may be identified isometrically with a subalgebra 
$B$ of $M_2(\ell^\infty(I))$ via the obvious map
$$
\left[ \begin{array}{ccl} \lambda I_H & y \\
0 & \mu I_H \end{array} \right]
\; \; \mapsto \; \; \left[ \begin{array}{ccl} \lambda 
1_I & i(y) \\ 0 & \mu 1_I \end{array} \right] \; \; \; .$$
Here $1_I$ is the identity of $\ell^\infty(I)$.  That this
is an isometry follows for example 
from \cite{FF} (Chapter IV \S2). It is
clear that this subalgebra $B$ is w*-closed, so that
$A$ is a dual Banach space.  Moreover the inclusion of $Y$
into  $A$ as the corner, is  a w*-homeomorphism.
That is, $Y$ is completely isometrically isomorphic to
a w*-closed subspace of $A$.
Therefore if $A$ with its given weak*-topology was 
w*-homeomorphic and completely isometric to
a $\sigma$-weakly closed subspace of some $B(H)$,
then so would $Y$ be,
which is false.

Finally, suppose that $Y$ is as in the example at the
end of \cite{EOR},
 an operator space which has a unique Banach space predual,
 but no operator space predual.  The $A$ constructed above
is a  dual Banach space.  If $A$ were a dual operator
space, there would exist by Theorem \ref{Blm}, a
completely isometric unital homomorphism $\pi$ of $A$ onto
a concrete dual operator algebra $B$.  However it is very easy
to check that the diagonal idempotents in $A$ force 
the Hilbert space on which $B$ acts, to split as a
direct sum $K \oplus N$, such that $B$ may be written as
$2 \times 2$ matrices exactly like our original
$A$.  Moreover, since there are  projections onto 
$K$ and $N$, the 1-2-corner of $B$ is  w*-closed in
$B(N,K)$.   Also, 
 $\pi$ restricted to the copy of $Y$, maps
$Y$ completely isometrically
onto this $1-2$ corner of $B$.
But this forces $Y$ to be a dual operator space.
\end{proof}

\vspace{5 mm}

Our theorem \ref{Blm}
 is thus best possible in the operator 
space category.   An obvious question which remains open
is the same question in the Banach category:
if $A$ is a unital operator algebra
which is a dual Banach space, then is
$A$ weak*-homeomorphic via an {\em isometric}
unital homomorphism to a $\sigma$-weakly closed subalgebra
of some $B(H)$?

Note that our characterization gives as a quick 
deduction, the fact due to Le Merdy \cite{LM3},
 and independently
Arias and Popescu \cite{AP}, that if 
$A$ is a concrete dual operator algebra (unital
or otherwise), and if $I$ is a weak*-closed 
two-sided ideal in $A$, then $A/I$ is
a dual operator algebra.  For if $A$ is
a weak*-closed subalgebra of $B(H)$, then 
by basic functional analysis so is
 $A + \Co I_H$, and $(A + \Co I_H)/I$ satisfies 
the hypotheses of our theorem.    This
gives the result.  However, this proof is not a
good proof for this result, since as Le Merdy has pointed out
to me privately, this $A/I$ result already follows
from the much simpler \ref{LM}.  Indeed one
need only note that the product on  $A/I$ is
automatically separately weak* continuous.  Le Merdy sees this 
by defining for each $b \in A$, a map $u : A \rightarrow A/I$ by 
$u(a) = [ab]$.  Clearly $u$ is  weak*-continuous.  Thus
its kernel is  weak*-closed, and we get an induced 
weak*-continuous map $A/(Ker \; u) \rightarrow A/I$.  
(Indeed the last sentence works for weak*-continuous 
continuous maps between any dual Banach spaces.)
The desired separate weak*-continuity follows immediately.
It is worth pointing out that this $A/I$ result whose proof we 
have just sketched is the major step in 
Le Merdy's earlier characterization of general dual operator algebras
(assuming separately weak*-continuous  product) in \cite{LM3}.
Thus Le Merdy's later characterization simplifies this
step.

It is interesting that the idea in the proof of
our theorem \ref{Blm} gives 
another proof of the following classical result.
Our argument has the advantage too of being applicable 
to dual function algebras (see comments after the 
proof).    

\begin{theorem} \label{nicep}  (Classical) Suppose that
$A$ is a $C(K)$ space,
or a commutative $C^*$-algebra,  with a predual.
Then $A$ 
 is w*-homeomorphically isometrically
*-isomorphic to a
commutative von Neumann algebra.
\end{theorem} 

\begin{proof}  
We may assume that $A = C(K)$ for compact $K$, 
by the usual argument
(if $A = C_0(K)$ for locally compact $K$ then
Krein-Milman implies the existence of extreme points.  Any
such extreme point, by a simple application of Urysohn's lemma 
has constant absolute value 1, so that $K$ is compact).  Replace
the argument above for Theorem
\ref{swc} by an almost identical argument in $B(A)$ instead
of $CB(A)$, and using the very easy \cite{HWW} Proposition I.3.9 instead
of our theorem \ref{BEZr}, to obtain that the multiplication
on $C(K)$ is separately w*-continuous.    Then
 we may use Le Merdy's
proof of \ref{LM} from \cite{LM4} to finish.
Alternatively,
there is a commutative proof  using measure theory
which avoids Le Merdy's argument.  We 
will omit this since this result is quite well known.
\end{proof}

\vspace{5 mm}

The argument above may be adapted, using the M-bounded 
characterization of multipliers \cite{HWW}, 
to any unital function algebra
(i.e. uniform algebra) with a predual, to show that any 
such algebra is w*- and isometrically isomorphic to
a commutative dual operator algebra.   Of course this fact 
also follows almost immediately from our theorem
\ref{Blm}.  This 
leads us to repeat from \cite{BLM}
 the following very interesting question
which seems to be open:  is a unital function algebra with
Banach space predual, w*- and isometrically isomorphic to
a w*-closed subalgebra of an $L^\infty$ space?  We know, by the
above, that it is a commutative dual operator algebra. 
If the  answer to this question is negative then  this
 would show that one must leave the category of function
algebras to study `dual function algebras'! 

Finally, 
we remark that Le Merdy's theorem \ref{LM} is true even without 
the presence of an identity (see \cite{LM3}); 
but we have no idea if the `separate
w*-continuity' hypothesis may be removed in the nonunital
case, except in the case that $A$ has a contractive approximate
identity.
Notice that if $A$ is an operator algebra 
with a contractive approximate
identity $\{ e_\alpha \}$,
and if $A$ has a predual and separately w*-continuous product,   
then if $e_\alpha \rightarrow e$ weak* in $A$, then $e$ is
clearly an identity in $A$ of norm 1.   However it seems
to be much more difficult to 
remove the `separately w*-continuity' hypothesis in the
last line.  We will do this by using an additional result from 
\cite{BEZ}.  In fact we have a more general result
(which does not assume `separate w*-continuity' of the 
product):

\begin{theorem} \label{rcai}  Let $A$ be an operator algebra
with a right contractive approximate
identity, and suppose that $A$ is a dual Banach space.
Then $A$ has a right identity of norm 1.
\end{theorem}

\begin{proof}  We know that $M_l(A)$ is a unital 
operator algebra, and it is easy to check that the 
natural map $L : A \rightarrow M_l(A)$ is a completely
isometric homomorphism.  Let $B = \{ T \in M_l(A) :
T L(a) = L(T(a)) \; \text{for all} \; \; a \in A \}$.
It is easy to see that $L(A) \subset B$, and that 
$B$ is a unital subalgebra of $M_l(A)$.   Clearly 
$B$ contains $L(A)$ as a left ideal.  Therefore by
\S 6 of \cite{BEZ}, $L(A)$ is a complete left M-ideal
in $B$.  By Theorem 3.10 (4) in \cite{BEZ}, 
$L(A)$ is a complete left M-summand of $B$.  Thus
$L(A) = Be$ for a projection $e \in B$.  Clearly
$e \in L(A)$, so that $L(A) = L(A)e$, from which
it follows that $e$ is a right identity for 
$L(A)$.
\end{proof}

\vspace{5 mm}

\section{Multiplier algebras of dual operator spaces}

We begin with a general functional analytic result and its 
operator space variant.  We will not need anything
above for its proof, except for \ref{KS} and \ref{mac}.

\begin{lemma} \label{semem} \begin{itemize}
 \item [(1).]  Let $X$ and $Y$ be  Banach spaces, with
$Y$ a dual Banach space, and let
$T : X \rightarrow Y$ 
be  a 1-1 linear 
map.  Then
the following are equivalent:
\begin{itemize}
\item [(i)]  $X$ is a dual Banach space and
$T$ is w*-continuous, and
\item [(ii)] $T(Ball(X))$ is w*-compact.
\end{itemize}

\item [(2).]  Let $X$ and $Y$ be operator spaces, with
$Y$ a dual operator space, and let
$T : X \rightarrow Y$  a 1-1 linear 
map
such that $T^{(n)}(Ball(M_n(X)))$ is w*-compact for every
positive integer $n$.  Then the 
predual of $X$ given in the proof of
(1), is an operator space predual of
$X$, and $T$ is w*-continuous.
\end{itemize}
\end{lemma}    

\begin{proof} (1) is undoubtedly well known,
we give a proof since we
have no reference for it.  Thanks go to 
C. Le Merdy for helping to simplify my original argument.
 That (i) implies (ii) is clear
(as indeed is the converse of statement (2), by the
way).  
Given (ii), 
we see by the Principle of Uniform 
Boundedness
that  $T(Ball(X))$ and therefore $T$ is bounded.   We
 may assume wlog that $T$ is a
contraction.  Suppose that $Z$ is the predual of $Y$,
and let 
$W = T^*(\hat{Z})$, a linear subspace of $X^*$.  The canonical 
map $j : X \rightarrow W^*$ is 1-1 and contractive and
$$\langle j(x) , T^*(\hat{z}) \rangle = T^*(\hat{z})(x)
= \hat{z}(T(x)) = T(x)(z) \; \; . $$ 
On the other hand, given $g \in Ball(W^*)$ let 
$\tilde{g}(z) = g(T^*(\hat{z}))$ for $z \in Z$, then
$\tilde{g} \in Ball(Z^*) = Ball(Y)$.  If $\tilde{g} = T(x)$
for an $x \in Ball(X)$ then we'd be done, for in this 
case its clear that 
$g = j(x)$.  So suppose, by way of contradiction
that $\tilde{g} \notin  T(Ball(X))$.
 By assumption $T(Ball(X))$ is w*-closed, so by 
the Hahn-Banach theorem,  
there exists $z \in Z$ such that 
$\tilde{g}(z) > 1$ and $|\langle T(x) , z \rangle|
\leq 1$ for all $x \in Ball(X)$.  The latter condition
implies that  $\Vert 
T^*(\hat{z}) \Vert \leq 1$, whereas the former condition
implies the contradictory assertion that
  $g(T^*(\hat{z})) > 1$.

That  $T$ is w*-continuous with respect to this predual 
of $X$ is now clear.

Now we prove (2).  We will use the fact (see \cite{LM1})
 that if an operator space $X$ is the dual of a
Banach space $W$, and if $W$ is equipped with its
natural matrix norms
 as a subspace of 
$X^*$ via the natural
inclusion, then $X$ is the dual operator space of
$W$ if and only if the unit ball of $M_n(X)$ is
$\sigma(X,W)$-closed.  To check that 
the latter condition holds, let
$x_\lambda = [x_{ij}^\lambda]$ be a  net in
$Ball(M_n(X))$, converging to $x =  [x_{ij}]$
in $M_n(X)$, the convergence being the
$\sigma(X,W)$-convergence.  That is,
$<x_{ij}^\lambda , w> \; \rightarrow \; <x_{ij},w>$
for all $i,j$ and $w \in W$.  Equivalently,
$< T(x_{ij}^\lambda) , z > \; \rightarrow \;
< T(x_{ij}) ,
z >$ for all $z \in Z$.  By Lemma \ref{mac},
the matrices $[T(x_{ij}^\lambda)]$ converge w*- to
$[T(x_{ij})]$ in $M_n(Y)$.  By hypothesis,
$[T(x_{ij})] \in T^{(n)}(Ball(M_n(X)))$, so that 
$x \in Ball(M_n(X))$, and we are done.
\end{proof}  

\vspace{5 mm}

In our situations below, we will know  in advance
that the operator $T$ is bounded, and so we
need only check w*-closedness instead of 
w*-compactness in the Lemma.
 
We will now show that the canonical inclusion 
$M_l(X) \rightarrow CB(X)$ satisfies the hypothesis
of (2) of the previous 
result\footnote{In fact a uniform version 
of this is true, namely that this map is what is known 
in the Banach space 
literature as a `semi-embedding', however we see no
use yet for this fact.}.  

\begin{corollary} \label{mis}  If $X$ is a dual operator
space then \begin{itemize}
\item [(1)]  $M_l(X)$ is a dual operator space.
\item [(2)]  A  bounded net $\{ a_i \}$
in  $M_l(X)$ converges
weak* to $a$  in $M_l(X)$ iff $a_i x \rightarrow
a x$ weak* in $X$ for all $x \in X$.
\item [(3)]  The product on  $M_l(X)$ is
separately weak*-continuous, and
$M_l(X)$ is a dual operator algebra.
\item [(4)]
If $M_l^w(X)$ is the set of $T \in M_l(X)$
which are weak*-continuous as maps on $X$,
then $M_l^w(X)$ is a unital norm closed subalgebra of
$M_l(X)$, and  the product on  $M_l^w(X)$ is
separately  weak*-continuous.
\end{itemize}
\end{corollary}

\begin{proof} (1):  Using \ref{semem}, we need to show 
that if $T_\lambda = [T_{ij}^\lambda]$ is
a net in $Ball(M_n(M_l(X)))$ converging 
w*- in $M_n(CB(X))$ to $T = [T_{ij}]$, then
 $T \in Ball(M_n(M_l(X)))$.  We will use the fact
from \cite{B} \S 4 that $M_n(M_l(X)) \cong
M_l(C_n(X))$, and we will test for
$T \in Ball(M_l(C_n(X)))$ using \ref{BEZr}.  So
we begin with two matrices $x = [x_{pq}(k)],
y = [y_{pq}(k)] \in M_m(C_n(X))$ (rows indexed
by $p$, columns indexed by $q,k$), and we need
to check that
$$\left| \left| \left[ \begin{array}{c}
\mu(T)^{(m)}(x) \\ y \end{array} \right]  \right|
\right| \; \; \leq \; \; \left| \left| \left[
 \begin{array}{c} x  \\ y \end{array} \right]  \right|
\right| $$
where $\mu : M_n(M_l(X)) \rightarrow M_l(C_n(X))$
is the canonical completely isometric identification.
However we do know that
$$\left| \left| \left[ \begin{array}{c}
\mu(T_\lambda)^{(m)}(x) \\ y \end{array} \right]  \right|
\right| \; \; \leq \; \; \left| \left| \left[
 \begin{array}{c} x  \\ y \end{array} \right]  \right|
\right| \; \; . $$
And of course $T_{ij}^\lambda(x) \rightarrow
T_{ij}(x)$ weak* in $X$ for all $i,j$ and
$x \in X$, by hypothesis.  The first matrix in the 
last displayed equation is an element of
$M_{2m,m}(C_n(X)) = M_{2mn,m}(X)$.  Now pick a norm 1
functional $G$ in the predual of the $M_{2mn,m}(X)$,
and apply it to this first matrix in displayed equation.
Using Lemma \ref{mac} as in the 
proof of \ref{wcl}
shows that
indeed $T \in Ball(M_l(C_n(X)))$.  Thus
we have now proved (1).
 
 Item (2) follows from (1)
and the definition of the weak*-topologies concerned
(see proof of \ref{semem}). 
Item (3) follows from Theorem \ref{Blm}.  

Finally, 
the  assertions of (4) are now fairly obvious.
If $T_n \in M_l^w(X)$, with
$T_n \rightarrow T$ in multiplier norm,
then $T_n \rightarrow T$ in cb-norm.  Since
the canonical image of $CB(X_*)$ in $CB(X)$
is norm closed, we see that $M_l^w(X)$ is
norm closed.
 \end{proof}

\vspace{5 mm}

Notice that the proof above provides, if
$X$ is a dual operator space such that 
the natural map $M_l(X) \rightarrow CB(X)$ is a 
complete isometry, a canonical
predual for $M_l(X)$ which is a natural
quotient of the operator space projective 
tensor product of $X$ and $X_*$.

The main part we wish to stress here is that
for any dual operator space
$X$, we have that $M_l(X)$ is a dual operator algebra,
that is it may be
w*-homeomorphically and completely isometrically
identified with a $\sigma$-weakly closed subalgebra
$B$ of some $B(K)$.

From this fact we can deduce another proof of the 
following result from \cite{BEZ}:

\begin{corollary} \label{enco}   If $X$ is
a dual operator space, then $A_l(X)$ is a W*-algebra.
\end{corollary}

\begin{proof}  From the above we now know that we can represent
$M_l(X)$ as a w*-closed subalgebra $B$ of some $B(K)$, with
$I_K \in B$.   Then the set of adjoints of operators in $B$ is
a w*-closed subalgebra of $B(K)$, and
$A_l(X) = \Delta(B) \subset B(K)$.   So $A_l(X)$ is
 w*-closed in $B(K)$, and consequently is a dual space.
\end{proof}

\vspace{5 mm}
                          
The last result was crucial for much of the theory in 
our project with Effros and Zarikian.   

We end this section with another application of 
`noncommutative M-structure'.  
We give 
a quick proof of an operator space variant on important
results in \cite{Zettl} and \cite{EOR}.  We recall that 
a TRO, or {\em ternary ring of operators},  bears the 
same relation to Hilbert $C^*$-modules, as 
closed *-subalgebras of $B(H)$ do to $C^*$-algebras.
Namely one may view TRO's as the `concrete' $C^*$-modules,
or view $C^*$-modules as `abstract' TRO's. 

\begin{corollary}  \label{wdu}  Suppose that
 $X$ is a TRO, ternary system, or right Hilbert
$C^*$-module, which is also a dual operator space.
Then \begin{itemize}
\item [(1)]  $X$ is a corner of a von Neumann algebra.
That is, $X \cong e R (1-e)$ completely isometrically and 
a weak*-homeomorphically,
for a von Neumann algebra $R$ and orthogonal projection $e \in R$.
\item [(2)]  If in addition $X$ is injective, then
$R$ may be taken to be injective too. 
\item [(3)]  Suppose that $X$ is full on the right over a
$C^*$-algebra $A$.  Let $M$ be the multiplier
algebra of $A$.   Then $M$, and
the algebra of adjointable operators on $X$, are
von Neumann algebras, and $X$ is a self-dual 
$C^*$-module over $M$ (and therefore also 
over $A$).  
\end{itemize}  
\end{corollary}

\proof
We prove (3) first.  We will use basic Hilbert $C^*$-module facts
\cite{La,Rieffelsu}.
We know from \cite{B}  
Appendix A.4 that 
for a right $C^*$-module 
$X$ we have that $A_\ell(X)$ is the
algebra of adjointable (in the usual $C^*$-module
sense) right module maps 
on $X$.  From \ref{enco}
$A_\ell(X)$  is a   von Neumann algebra.
By Morita symmetry, and using Kasparov's result asserting that 
the algebra of adjointable maps is the multiplier algebra of the 
algebra of `compact' adjointable maps, we have that 
the multiplier algebra $M$ of $A$
 is isomorphic to
the algebra $A_r(X)$ of adjointable operators acting on the 
right of $X$.  By
symmetry this is also
a  von Neumann algebra.  By \ref{dn} in the next section, 
the product map $A_\ell(X) \times X \times M \rightarrow X$
is separately weak*-continuous.
  We will now check that the inner product
$X \times  X \rightarrow M$ is separately weak*-continuous.
To this end suppose that we have
a bounded net $y_\lambda \rightarrow y$ weak* in $X$.  Then
letting $T$ be the rank 1 adjointable operator $w \otimes x
$ for $x, w \in X$, we have by the above that
$T y_\lambda \rightarrow Ty$ weak*, so that
$w \langle x ,  y_\lambda \rangle \rightarrow w \langle x , y
\rangle$  weak*.
The net $\{  \langle x ,  y_\lambda \rangle \}$ is bounded, and if it has
a weak*-convergent subnet $\{  \langle x ,  y_{\lambda_\mu} \rangle \}$,
converging to $b \in M$ say, then by the first part
$w \langle x ,  y_{\lambda_\mu} \rangle  \rightarrow wb$ weak*.  Hence
$wb =  w \langle x , y \rangle$.  Since this is true for all $w \in X$,
it follows that $b = \langle x , y \rangle$.  Hence
the net $ \langle x ,  y_\lambda \rangle \rightarrow \langle x , y
\rangle$  weak*.
 
We now can see that $X$ is a self-dual right 
$M$-module by the following trick.
Suppose that $f : X \rightarrow M$
is a bounded $M$-module map.
It is well known that one may choose
a contractive approximate identity $\{ T_\alpha \}$ for
the `imprimitivity $C^*$-algebra' $K_A(X)$
with terms of the form $\sum_{k = 1}^n x_k \otimes x_k$ for
$x_k \in X$.   Note that $f(T_\alpha(x)) =
\sum_{k = 1}^n f(x_k) \langle x_k , x \rangle =
\langle \sum_{k = 1}^n  x_k f(x_k)^* , x \rangle$.  The element
$\sum_{k = 1}^n  x_k f(x_k)^*$ depends on $\alpha$, let us
name it $w_\alpha$.  Then $\{ w_\alpha \}$ is a bounded net in
$X$, so has a subnet converging weak* to a $w \in X$.  By replacing the
net with the subnet, and using the
first part, $\langle w_\alpha , x \rangle \rightarrow \langle
w , x \rangle$.  Since
$f(T_\alpha(x)) \rightarrow f(x)$  in norm, we obtain
$f(x) =  \langle w , x \rangle$
for all $x \in X$.   So $X$ is  self-dual over $M$.

(1) follows immediately from
(3) as is well known
\cite{Paschke,Rieffel}.  
For example this is easily seen
by considering the linking
von Neumann algebra,
which may be taken to be
$B_M(X \oplus_c M)$, where $X \oplus_c M$ is is the
$C^*$-module direct sum).  Indeed the $C^*$-module direct sum
$X \oplus_c M$ is clearly selfdual over $M$, and hence
by \cite{Paschke}, $B_B(X \oplus_c M)$ is a von Neumann algebra.

(2): If $X$ is injective in addition,
then we know from \cite{BPnew} (see the introduction and
Corollary 1.8 there) that  (using the notations of that paper) 
$I_{22}$ is the multiplier algebra of $X^*X = A$,
so that $I_{22} = M$.
Similarly, the linking von Neumann
of the last paragraph, 
is what we called $I(S(X))$ in
\cite{BPnew} 
 and is therefore injective.
So $X$ is a corner of an injective von Neumann algebra.
\endproof%

There are some similar arguments to the above in the
thesis of J. Schweizer \cite{Schw}.  We also 
observe that another proof of 
the last corollary proceeds on the following lines.
If $X$ is a dual operator space then
 the natural projection $P :
X^{**} \rightarrow X$ is completely contractive.
One then sees from results in \cite{Yo} that 
$P$ is a `conditional expectation'.  This allows
one to circumvent
a large part of the argument in \cite{EOR} (which 
is exactly what one would expect, given our stronger 
hypothesis).

\section{Characterizations of dual actions on an
operator space}

We will use the following simple result from \cite{BEZ}.  For 
completeness we give a short proof.  

\begin{proposition}
\label{tem}
If $X$ is a dual operator space then any
$T \in A_l(X)$ is weak*-continuous.
\end{proposition}

\begin{proof}   We may assume that $0 \leq T \leq 1$
in the $C^*$-algebra $A_\ell(X)$.  By basic operator theory,
any such $T$ is the 1-1-corner of an orthogonal
  projection $P \in M_2(A_\ell(X))$.  However
$M_2(A_\ell(X)) \cong A_\ell(C_2(X))$.  Moreover the 
left adjointable projections are exactly the 
complete left $M$-projections, which are shown 
in  \cite{BEZ} to be weak* continuous
(the argument for this is quite simple and
 very similar to the classical one).
Thus $P$ is weak*-continuous.
We have $T = \pi_1 P j_1$,
where $\pi_1 : C_2(X) \rightarrow X$
and $j_1: X \rightarrow  C_2(X)$ are the canonical maps.
Since these
latter maps are clearly weak*-continuous, so is $T$.      
\end{proof}

\begin{corollary}  \label{dn} Any dual operator space $X$ is a
normal dual $A_l(X)-A_r(X)$-bimodule.  Also the left
module action of $M_l(X)$ or $A_l(X)$
on $X$ is weak*-continuous in the
first variable.
\end{corollary}

The proof of the last corollary is elementary from
\ref{tem} and the definition of the
w*-topology on $A_l(X)$ from \S 3 (i.e. Corollary \ref{mis}
(2)).

Unfortunately, at this point we do not know whether 
the action of $M_l(X)$ on $X$ is weak*-continuous in
second variable, or equivalently whether left multipliers 
on $X$ are automatically  weak*-continuous on $X$.  This
makes it necessary to make a definition: 
we say that a dual operator space is {\em left
normal} if every $T \in M_l(X)$ is weak*-continuous
on $X$.   Unital dual operator algebras, dual
TRO's and reflexive spaces
are all left normal, for example.            

\begin{theorem} \label{meg}  Let  $m : M \times X 
\rightarrow X$ be an
operator module action of a
W*-algebra on a dual operator space, or more generally an
oplication of a dual operator system on a dual operator space.
If $m$ is weak*-continuous in the first variable,
then it is separately  weak*-continuous.  Moreover, in this
case there exist
 Hilbert spaces $H, K$, 
 a unital w*-continuous completely positive 
map   $\pi : M \rightarrow B(K)$, and a
weak*-continuous complete isometry $\Phi : X \rightarrow
B(H,K)$, such that $\Phi(m(b,x)) = \pi(b) \Phi(x)$
for all $x \in X, b \in M$.
Moreover $H, K, \Phi$ can be chosen to
only depend on $X$, and not on $M$ or
the particular action.  If $M$ is a W*-algebra,
then $\pi$ is a
*-homomorphism if and only if $m$ is a module action.
\end{theorem}           

\begin{proof}  First apply Theorem \ref{ern} to 
the $A_l(X)$ action on $X$ (using Corollary \ref{dn}).
This provides the Hilbert spaces  $H, K$ and the map $\Phi$.
We also obtain a normal *-homomorphism $\rho
 : A_l(X) \rightarrow B(K)$, with $\rho(T) \Phi(x) = 
\Phi(Tx)$ for all $x \in X, T \in A_l(X)$.  Next apply the
oplication theorem,
to get a unique
 unital completely positive $\psi : M \rightarrow A_l(X)$
such that $\psi(b) x = m(b,x)$ for all $x \in X, b \in M$. 
The last assertion in the 
statement of our theorem follows also from the 
oplication theorem.
It is elementary to check from the definition of the 
w*-topology on $A_l(X)$ from \S 3, the fact that 
$m$  is weak*-continuous in the first variable,
and (2) of the Krein-Smulian theorem, that $\psi$ is
weak*-continuous.  Let $\pi = \rho \circ \psi$ . 
Finally, if $b \in M$, and 
$x_i \rightarrow x$ w*- in $X$, then
$\pi(b) \Phi(x_i)  \rightarrow \pi(b)
\Phi(x)$ w*- in $B(H,K)$.  That is,
$m(b,x_i) \rightarrow m(b,x)$ (since $\Phi$ is a 
w*-w*-homeomorphism).
Thus $m$ is weak*-continuous in the second 
variable.
\end{proof}
 
\vspace{5 mm}

We should also remark on the necessity of the hypothesis
that $m$ be weak*-continuous in the
first variable, in the above theorem.
To see this consider the natural action of $B(H)$ on $B(H)^{**}$.
This is an operator module action which is clearly not
weak*-continuous in the
first variable.

We now turn away from the case where $M$ is a W*-algebra or 
dual operator system, and we will replace it by a dual operator
algebra, or general dual operator space.  Thus we 
consider oplications $A \times X \rightarrow X$, where 
$A, X$ are dual operator spaces.  One immediately has:

\begin{theorem} \label{ndo}
Suppose that $m : A \times X \rightarrow X$ is an
oplication of a
dual operator space $A$ on a dual operator space
$X$, which is weak*-continuous in the
first variable.  Then there exists a unique
weak*-continuous unital completely
contractive $\pi : A  \rightarrow M_l(X)$ such that
$\pi(a) x = m(a,x)$ for all $a \in A, x \in X$.
Also:
\begin{itemize}
\item [(1.)]  If $A$ is an algebra then
$\pi$ is a homomorphism iff $m$ is a module
action.
\item [(2.)]  If $m$ is separately w*-continuous
then
the range of $\pi$ is contained in $M_l^w(X)$.
\end{itemize}
\end{theorem}

\begin{proof}  This follows immediately from the  
oplication theorem 
which provides us with a
unital $\pi : A \rightarrow M_l(X)$ such that
$\pi(a) x = m(a,x)$ for all $a \in A, x \in X$.
As in \ref{meg}  $\pi$ is weak*-continuous.   
A similar argument proves (2.).
\end{proof}                            

\vspace{5 mm}

The following characterizes dual oplications as corresponding
to appropriate concrete
weak*-closed spaces of operators on Hilbert space.  It is
a generalization of Theorem 3.1 in \cite{BLM} (which
gave a similar characterization for dual operator modules
over dual operator algebras):    

\begin{theorem} \label{meg2} Suppose that
 $m : A \times X \rightarrow X$ is an oplication of a
dual operator space $A$ on a  dual operator space
$X$, with $e$ the `identity' of $A$.  Suppose that 
either $m$ is separately weak*-continuous, or that 
$X$ is left normal and $m$ is  weak*-continuous
in the first variable.
Then there exist Hilbert spaces $H$ and $K$
a weak*-homeomorphic
completely isometric embedding  $\Phi : X \rightarrow B(H,K)$,
and a weak*-continuous completely contractive
linear map $\theta : A \rightarrow B(K)$, such that
$\theta(e) = I_K$, and
$\Phi(m(a,x)) = \theta(a) \Phi(x)$ for all $a \in A, x \in X$ .
If $A$ is also an algebra, and $m$ is a
module action then
$\theta$ is a homomorphism.  
Finally,
the map $\Phi$ and spaces $H$ and $K$ may be chosen to only
depend on $X$, i.e. they are independent of $A$ and $m$.
             \end{theorem}

\begin{proof}  By \ref{ndo}
there
is a unique w*-continuous
unital complete contraction $\rho : A  \rightarrow
M_l^w(X)$ such that $m(a,x) = \rho(a) x$.  If $X$ is left
normal then the proof proceeds exactly as in \ref{meg},
but using Theorem 3.1 in \cite{BLM} instead of
\ref{ern}, 
to obtain the result.  So we assume henceforth that 
$m$ is separately weak*-continuous.

Let $B = M_l^w(X)$, a closed unital subalgebra of the
dual operator algebra $M_l(X)$.   On $B$ we consider
the relative
weak*-topology inherited from $M_l(X)$.
Consider the natural operator module action
$n : B \times X \rightarrow X$.  Then $n$ is
separately weak*-continuous.
If one follows through
the proof of Theorem 3.1 in \cite{BLM} carefully one finds
that, even without assuming that  $B$ is a dual space,
we have a weak*-homeomorphic
completely isometric embedding  $\Phi : X \rightarrow B(H,K)$,
and a completely contractive unital
homomorphism $\pi : B \rightarrow B(K)$, such that
$\Phi(n(b,x)) = \pi(b) \Phi(x)$ for all $b \in B, x \in X$ .
Moreover for any bounded net $b_\lambda \rightarrow b$   
weak*- in $B$, we have $\pi(b_\lambda) \rightarrow \pi(b)$
weak*- in $B(K)$.

Let $\theta = \pi \circ \rho$.   Then $\theta$ restricted
to $Ball(A)$ is  w*-continuous, so  $\theta$ is  w*-continuous.
The rest is clear.
 \end{proof}

\vspace{5 mm}

There is a bimodule version of the last theorem which we 
leave to the reader.   A quite interesting feature of 
these theorems, and \ref{meg}, is that it provides 
for any dual operator space $X$, 
a completely isometric w*-w*-homeomorphic representation
of $X$ inside $B(H,K)$ for some Hilbert spaces, and 
dual operator algebras inside $B(K)$ and $B(H)$, which are 
universal in that every `normal dual bimodule' action
on $X$ has to `factor through' this one.

Theorem \ref{meg2}
is not quite a complete analogue of 
Theorem \ref{meg}, in that we were not able to answer the
question: does
w*-continuity in the first variable imply w*-continuity in the 
second?  But it turns out that this is exactly the same as
the question of the existence of a non-weak*-continuous
left multiplier on $X$:

\begin{corollary}  Let $X$ be a dual operator
space.  There exists a left
operator module action of
a dual operator algebra on $X$ 
which is 
weak*-continuous in the first variable but 
not the second, if and only
if there
exist left multipliers which are not weak*-continuous on $X$.
We are not sure if this ever happens.
\end{corollary}

\begin{proof}  If $X$ is not left normal, then the
action of $M_l(X)$ on $X$  is only
weak*-continuous in the left variable.  Conversely,
if $X$ is  left normal, and we are given a left module action 
of a dual operator algebra on $X$ which is 
weak*-continuous in the left variable, then it is
weak*-continuous in the second variable exactly as in \ref{meg},
but using \cite{BLM} Theorem 3.1 instead of
\ref{ern}.  \end{proof}

\vspace{5 mm}
                   
Next we make a remark on unital operator spaces which 
are also dual operator spaces.  In this case, as explained in the
introduction, the map $M_l(X) \rightarrow X$ given by
$T \mapsto T(1)$ is a unital complete isometry.  However it
is clearly weak*-continuous with the topology on
$M_l(X)$ from \S 3, so that by
\ref{KS}, it is a weak*-homeomorphism onto its 
w*-closed range $W$.  Thus $X$ contains a 
w*-closed unital subspace which is also a dual operator
algebra.  It is clear that
any $T \in M_l(X)$ has the property that $T(W) \subset W$,
so that $T$ restricted to $W$ is a multiplier of 
$W$ by Corollary \ref{her}.  Consequently $T$ is 
weak*-continuous on $W$, since $W$ is a dual operator 
algebra.
 
Finally we consider the `central 2-sided' version of the 
theory of this section.   We will assume familiarity 
with some basic ideas from \cite{B}.    
 Recalling that $A_l(X)$ corresponds
to left operator module actions of a $C^*$-algebra on 
$X$, and that $A_r(X)$ corresponds to right actions, one 
is led to consider the subset $A_l(X) \cap A_r(X)$ inside
$CB(X)$.  This space is a commutative $C^*$-algebra which
we have thoroughly studied elsewhere, and for which
we have 
quite a large number
of characterizations.  We call it the {\em operator space
centralizer algebra} $Z(X)$, not to be confused with the
Banach space centralizer algebra $Cent(X)$ (see \cite{HWW}).
Nonetheless, $Z(X)$ may be developed entirely analogously to 
the classical centralizer theory \cite{AE,HWW}.  The ensuing theory
might be called  `central complete M-structure', as opposed to 
`complete right M-structure'.  
`Central complete M-structure'  is much
less interesting in some ways than the  
`1-sided' theory precisely because it is so close to the
classical, commutative theory surveyed 
in \cite{HWW}.  We will present
these details elsewhere.   Here we will simply discuss
the `central versions'
 of the characterization theorems of the type
above.  

We first note that in the development of 
$A_r(X)$ as opposed to $A_l(X)$ there is a slight twist,
the subset $A_r(X)$ of right adjointable maps on $X$ should be 
given the opposite of the usual multiplication of $CB(X)$
if we want it to 
correspond *-isomorphically to a subalgebra of ${\mathcal
F}(X)$ in the language of \cite{B} (or $I_{22}(X)$
in the language of \cite{BPnew}).  This will not really
be an issue for us though.  It follows easily by looking
at these latter algebras, that $S T = T S$ if $S \in A_l(X),
T \in A_r(X)$.   Thus the subalgebra
$Z(X) = A_l(X) \cap A_r(X)$ of $A_l(X)$ is a $C^*$-algebra
which is also commutative.  It also follows from 
\ref{enco}   now that  
if $X$ is a dual operator space then
 $Z(X)$ is a commutative W$^*-$algebra.

We will say that a left operator module $X$ over a 
$C^*$-algebra $A$ is a {\em central operator module} if
the map $X \otimes_h A \rightarrow X$ given by
$(x,a) \mapsto ax$, is also completely contractive.  Or, in
other words, this latter map is a right oplication.  
If 
$A$ is commutative, this is just saying that $X$ is an
operator $A-A$-bimodule with respect to this action
$a x = x a$.   It should be quite clear now that any 
operator space $X$ is a central operator
$Z(X)$-module.    Conversely, it is easy to see from 
the oplication theorem that for  
any central operator $A$-module $X$ there exists a
*-homomorphism $\pi : A \rightarrow Z(X)$ such that 
$\pi(a) x = a x$ for all $a \in A, x \in X$.  If 
$A$ is unital then so is $\pi$.    

There is a similar definition for `central oplications',
and an analogous argument showing that 
every  central oplication can be written
as in the last equation, for a completely positive $\pi$.
 
It is clear that one way to produce central operator 
modules is to consider subalgebras of the commutant $X'$
of a concrete operator space $X$.  More particularly,
if $X \subset B(H)$, and if
$B \subset X' \subset B(H)$ is such that 
$BX \subset X$, then $X$ has a central $B$-action.  It is 
fairly 
obvious how to show that these are essentially the only 
ones; namely by first applying Theorem 2.1 of 
\cite{ERbimod} to the $Z(X)$-bimodule action on
$X$, and then composing the *-representation obtained 
there with the $\pi$ two paragraphs above.   Here, in more detail,
 is the weak*-version of such a result:

\begin{theorem} \label{cen}  Let $B$ be a W$^*-$algebra (resp.
dual operator system), and let $X$ be a dual operator space,
and let $m : B \times X \rightarrow X$ be a central
 operator module action (resp. central oplication) which is
w*-continuous in the $B$-variable.  Then there exists a
Hilbert space $H$, a
 w*-homeomorphic complete isometry $\Phi : X 
\rightarrow B(H)$, and a unital weak*-continuous
*-homomorphism (resp. completely positive map)
$\pi : B \rightarrow \Phi(X)' \subset B(H)$ such that
$\pi(b) \Phi(x) = \Phi(m(b,x))$  for all $b \in B, x \in X$. 
As above, $H, \Phi$ may be chosen independently of
$B, m$.
\end{theorem}

\begin{proof}  The idea is the same as for the proof
of \ref{meg}.  We first observe that the $Z(X)$ action on
$X$ (which clearly makes $X$ a  
normal dual $Z(X)$-bimodule by \ref{dn})
  has such a representation.  To see this use  
\cite{ERbimod} Theorem 4.1 
that there exist such $H, \Phi, \pi$ as above, 
$\pi : Z(X) \rightarrow B(H)$, such that 
$\pi(z_1) \Phi(x) \pi(z_2) = \Phi(z_1 x z_2)$ for all
$z_1, z_2 \in Z(X)$, $x \in X$.  
Then, since 
$\pi(z) \Phi(x) = \Phi(x) \pi(z)$ we have that 
$\pi$ maps into $\Phi(X)'$.  Now proceed as in 
\ref{meg}.   We leave the omitted details to the reader.
\end{proof}  

\vspace{5 mm}

We now classify the singly generated central modules:

\begin{theorem} \label{sg}  Let $X$ be a central operator
$A$-module over a $C^*$-algebra $A$.  Suppose that there exists
an $x_0 \in X$ such that the closure of $A x_0$ is norm dense 
in $X$ (or weak*-dense in $X$, if 
$X$ is a dual operator space).  Then $X$ is a $MIN$ space, so that 
$X$ is a function $A$-module in the sense of \cite{BLM}.  That 
is, there is a compact space $K$, and a
completely isometric 
linear $ \Phi : X \rightarrow C(K)$,
and a *-homomorphism $\pi : A \rightarrow C(K)$, 
such that $\Phi(ax) = \pi(a) \Phi(x)$ for all $a \in A, x \in X$. 
Conversely, any function $A$-module is a central operator
module.  If $X$ is an algebraically singly generated
(i.e. $A x_0 = X$)  central operator module over a $C^*$-algebra
 then in fact 
$X$ is a $C(K)$ space (linearly completely isometrically).
\end{theorem} 

\begin{proof}   Suppose that $A x_0$ is norm dense in $X$.
By the result three paragraphs above \ref{cen}, 
it follows that 
$Z(X) x_0$ is norm dense in $X$.   It is easy to see 
that $Z(X)$ commutes
with $J(x_0) J(x_0)^*$ inside the multiplier 
algebra of the $C^*$-algebra
we called ${\mathcal E}(X)$ in \cite{B}.  We refer to that
paper, particular \S 6 there,
 for background on what follows.
Since $Z(X) J(x_0)$  
is norm dense in $J(X)$, it follows 
from the definition of
${\mathcal E}(X)$, that the latter will be a
commutative $C^*$-algebra.  Note also, that since 
$X$ is also singly generated as a right $Z(X)$ module,
if
$E = (J(x_0) J(x_0)^*)^{\frac{1}{2}}$, then by
\cite{B} Lemma 6.2, $E$ is strictly positive in
${\mathcal E}(X)$.  Hence
 the set $E {\mathcal E}(X)$ is dense in ${\mathcal E}(X)$,
by Stone-Weierstrass.  
By \cite{La} Lemma 4.4, we may write 
$J(x_0) = E^{\frac{1}{2}} w$, for some
$w \in {\mathcal T}(X)$.
Since ${\mathcal E}(X) J(X)$ is dense in ${\mathcal T}(X)$
(see the beginning of \S 4 in
\cite{B}), the map $Ez \mapsto z E^{\frac{1}{2}} w$ from
$E {\mathcal E}(X) \rightarrow {\mathcal T}(X)$,
is  an isometric
${\mathcal E}(X)$-module map onto a dense
subset of ${\mathcal T}(X)$.  Hence it extends
to a completely isometric surjective
${\mathcal E}(X)$-module map
$T$ say,
from ${\mathcal E}(X)  \rightarrow {\mathcal T}(X)$.
Thus ${\mathcal T}(X)$ and therefore also
$X$, are Banach spaces with the MIN structure.
In this case we are in the setting of \cite{BLM}
and \cite{B} \S 2 and 3.  That is, $X$ is a topologically
singly generated function module over $Z(X)$.  We see from
Corollary 3.6 in \cite{B} that $X$ is algebraically
singly generated iff $X$ is e.n.v. as a Banach space,
and in this case $X$ is actually a commutative unital
$C^*$-algebra.

Finally, if $X$ is a dual space and $A x_0$ is w*-dense,
then the above shows that the norm closure of 
$Z(X) x_0$ is a  MIN space.  However the
weak*-closure of a MIN space is also  a MIN space
(we leave this fact as an exercise), so that
  $X$ is a MIN space.
\end{proof}

\vspace{5 mm}

Thus the topologically singly generated central operator 
modules over a  $C^*$-algebra $A$ are exactly the  
topologically singly generated function modules over $A$.  

Notice that every element $x_0$ in an operator space $X$,
is contained in a subspace of
$X$ which is a topologically singly
generated central operator $Z(X)$-module, namely the closure
of $Z(X) x_0$.   We know this is a $MIN$ space by the last 
result, but no doubt it will in most cases be 1-dimensional.   

We also remark that the result above that an algebraically
singly generated central operator module (or
equivalently, function module) over a $C^*$-algebra 
is linearly isometric to
 a commutative $C^*$-algebra, is false with 
`algebraically' replaced by `topologically singly 
generated'.    Simple examples abound (see remark after
3.6 in \cite{B}).

\vspace{5 mm}

One final remark:  there are a host of well known
characterization theorems which are appropriate
for the
`completely isomorphic' case of the theory, as opposed to the
`completely isometric' type results above
(see \cite{BLM,Bnat,Pis,LM3} for
example).   Of course our techniques here have no hope of working
in the `completely isomorphic' case.

\end{document}